\documentclass[11pt,thmsa]{article}
\usepackage{amsfonts}
\usepackage{amssymb}
\newtheorem{teo}{Theorem}[section]

\newtheorem{lema}[teo]{Lemma}
\newtheorem{cor}[teo]{Corollary}
\newtheorem{prop}[teo]{Proposition}
\newtheorem{defin}[teo]{Def\mbox{}inition}
\newtheorem{obs2}[teo]{Remark}

\newtheorem{tea}{Theorem}[subsection]

\newtheorem{no2}[teo]{Note}

\newtheorem{no3}[tea]{Note}

\newcommand{\Gal}{{\rm Gal}}
\newcommand{\Frob}{{\rm Frob }}
\newcommand{\trace}{{\rm trace}}
\newcommand{\mod}{{\rm mod}}
\newcommand{\numerator}{{\rm numerator}}

\newcommand{\disc}{{\rm disc}}

\newcommand{\SL}{{\rm SL}}

\newcommand{\GL}{{\rm GL}}
\newcommand{\Norm}{{\rm Norm}}

\newcommand{\PSp}{{\rm PSp}}
\newcommand{\PGSp}{{\rm PGSp}}

\newcommand{\GSp}{{\rm GSp}}
\newcommand{\Sp}{{\rm Sp}}

\newcommand{\End}{{\rm End}}
\newcommand{\diag}{{\rm diag}}

\title{On the images of the Galois representations attached to genus $2$
Siegel modular forms}
\author{Luis V. Dieulefait  \thanks{supported by TMR - Marie Curie
Fellowship ERBFMBICT983234} }

\begin{document}

\maketitle
\begin{abstract}
We address the problem of the determination of
the images of the Galois representations attached to genus $2$ Siegel
cusp forms of level $1$  having multiplicity one. These representations are
 symplectic. We prove
 that the images are as large as possible for
almost every prime, if the Siegel cusp form is not a Maass spezialform and
verifies two  easy to check conditions.
\end{abstract}

\section{Introduction}

In this article we apply the same techniques as in [5]
to study the images of the  Galois representations
attached to genus $2$ cuspidal Siegel modular forms on
$\Sp(4, \mathbb{Z})$ of even weight $k$ that are Hecke eigenforms.
 These four dimensional Galois representations have
been constructed by Taylor and Weissauer (see [18], [22]).\\
We will restrict ourselves to the case where the automorphic representation
corresponding to the Siegel modular form has multiplicity one. In this
case the Galois representations are symplectic (cf. [22], [3]).

\medskip
\noindent We will consider explicit examples of such Siegel modular forms whose
first Fourier coefficients have been computed by Skoruppa (see [15], [16])
and we will show that the images are ``as large as possible"  (a certain
symplectic group) for almost
every prime. The examples of Skoruppa do not correspond to Maass spezialformen
 (Siegel
modular forms in the image of
the Saito-Kurokawa lift from $S_{2k-2}(\SL_2(\mathbb{Z}) )$),
 and in addition  Eisenstein-Klingen series
are excluded. For such cases it is known that the image is not
maximal for any prime.

\medskip
\noindent Our main tools will be
the classification of maximal subgroups of $\PSp(4,
\mathbb{F}_q)$ (see [10]); the description of the
 image of the inertia subgroup $I_\ell$ due to the work of
 Faltings  (see [6], [2] and [20])  on Hodge-Tate
 decompositions  via an application of Fontaine-Laffaille theory;
  and Serre's
 conjecture $ (3.2.4 _?)$ (cf. [14]) for odd irreducible two-dimensional Galois
 representations.

\medskip
\noindent Our main result shows, under the assumption of Serre's conjecture
 $(3.2.4_?)$, that
 given a genus $2$  Siegel cusp form for the full Siegel modular group
 of even weight, if it has multiplicity one,  is not a Maass spezialform
  and verifies
 an additional condition (similar to the condition
 ``without inner twists" imposed in [4] to classical modular
 forms), then the images of the corresponding Galois
representations are ``as large as possible" for almost every prime.\\
In fact, a stronger statement is proved in the last section using the theory
 of pseudo-representations: the same
 result holds independently of Serre's conjecture if
we further impose the condition that for some prime $p$:
 $ \mathbb{Q}(d_p) = E $ and $ \sqrt{d_p} \notin
  E $ where $E$ is the number field generated by the eigenvalues of the
  Siegel cusp form and $d_p$ is defined in section 4.8.\\
This condition is verified in the examples. Moreover,
one can easily prove assuming Serre's conjecture
that this condition will always be satisfied.\\
I want to thank N£ria Vila for
many helpful comments and suggestions.

\section{The setup}

Let $f$ be a cuspidal weight $k$ Siegel modular form for $\Sp(4, \mathbb{Z})$
which is a Hecke eigenform and denote $a_n$ its eigenvalues.\\
Let
$$ Z_f(s):= \zeta(2s-2k+4) \sum_{n=1}^{ \infty } \frac { a_n} {n^s} $$
be the spinor zeta function. Then $Z_f$ has an Euler product of the
form
$$ Z_f(s) = \prod Q_p (p^{-s})^{-1} ,$$
where $Q_p$ is the polynomial:
$$ Q_p (x) = 1 - a_p x + (a_p^2 - a_{p^2} - p^{2k-4})x^2 - a_p p^{2k-3} x^3
+ p^{4k-6} x^4 .$$
Then we have the following
  result of Taylor [18], completed by Weissauer [22] (see also [3]):
\begin{teo}
\label{teo:Weissauer} Let $f$ be a Siegel modular form of even weight $k$
on the full Siegel modular group $\Sp(4, \mathbb{Z})$ which is a cusp
form and a simultaneous eigenform for all Hecke operators $T(n)$.
Furthermore, assume that the automorphic representation corresponding
to $f$ has mutiplicity one. Let $E$
be the number field generated by the eigenvalues $a_n$. For any prime
number $\ell$ and any extension $\lambda$ of $\ell$ to $E$ there exists
a continuous Galois representation
$$ \rho_{f,\lambda} : G_{\mathbb{Q}} \rightarrow \GSp(4, \overline{E}_\lambda)$$
such that the following holds: the representation $\rho_{f,\lambda}$ is
unramified outside $\ell$ and
$$ \det( Id - x \rho_{f,\lambda}(\Frob \; p)) = Q_p (x)$$
for every $p \neq \ell$. If $ \rho_{f, \lambda} $ is absolutely
irreducible, then it is defined over $E_\lambda$.\\
\end{teo}
For such an $f$ and under the hypothesis of
absolute irreducibility, the maximal possible image for
$ \rho_{f, \lambda}$ is the group
$$ A_{\lambda}^k = \{ g \in \GSp( 4, \mathcal{O}_{E_\lambda}) :
\det(g) \in ( \mathbb{Q}_\ell^*)^{4k-6} \} ,$$
where  $ \mathcal{O}_{E_\lambda} $ is the ring of integers of
$E_\lambda$.\\
If we consider the residual $\mod \; \lambda$ Galois representation
$\bar{\rho}_{f,\lambda}$ (see the remark below), then its  image
is contained in
$$  \{ g \in \GSp( 4, \mathbb{F}_{\lambda}) :
\det(g) \in ( \mathbb{F}_\ell^*)^{4k-6} \}  .$$
The image of   its projective image $\mathbb{P}
( \bar{\rho}_{f,\lambda}  )$ is contained
 in $ \PSp(4, \mathbb{F}_{\lambda})$, if the residual degree of $\lambda$
is even, and in $ \PGSp(4, \mathbb{F}_{\lambda})$ if it is odd (see
[3]).

\medskip
\noindent {\bf Remark:} In these assertions  the assumption of absolute irreducibility is
not needed.
 This follows from the fact that $E$ is the field generated by
the coefficients of the characteristic polynomials of the images of the
Frobenius elements and from a lemma on the reduction of $\ell$-adic representations
proved in [3]. This lemma tells us that the $ \mod \;
\lambda$ residual representation obtained by reducing the coefficients of
all the characteristic polynomials of the images of Frobenius elements
  ``agrees" with the
reduction of our $\lambda$-adic Galois representation that may not be defined
over $E_\lambda$; actually it is only thanks to
 this lemma that we can define  the residual $\mod \; \lambda$
representation  in general.\\
Once we have determined that for a prime $\ell$ the image of
$\bar{\rho}_{f,\lambda}$ is as large as possible, an application of a
lemma of Serre (cf. [13] and [3]) implies that the image of
$ \rho_{f,\lambda}$ is  $  A_{\lambda}^k  $.\\

\section{The tools}

\subsection{Maximal subgroups of $ \PGSp(4, \mathbb{F}_q) $ }

In [10], Mitchell gives the following classification of maximal
proper subgroups $G$ of $\PSp(4, \mathbb{F}_q)$ where $q=p^r $, p odd prime, as
groups of transformations of the projective space having an invariant
linear complex:\\
$ $\\
1) a group having an invariant point and plane\\
2) a group having an invariant parabolic congruence\\
3) a group having an invariant hyperbolic congruence\\
4) a group having an invariant elliptic congruence\\
5) a group having an invariant quadric\\
6) a group having an invariant twisted cubic\\
7) a group $G$ containing a normal elementary abelian subgroup $E$ of
order $16$, with: $ G/E \cong A_5 \; \mbox{or} \;  S_5$ \\
8) a group $G$ isomorphic to $A_6 , S_6$ or $ A_7$\\
9) a group conjugated under $\PSp(4, \mathbb{F}_{p^r})$ with
$\PSp(4, \mathbb{F}_{p^k})$, where $r/k$ is an odd prime\\
10) a group conjugated under $\PGSp(4, \mathbb{F}_{p^r})$ with
$\PGSp(4, \mathbb{F}_{p^k})$, where $r$ is even and $r/k = 2$\\
$ $\\
Cases 7) and 8) only occur if $r=1$.\\
(For the relevant definitions see [8], see also [1] and [11]
for cases 7 and 8).\\
\noindent
From this we obtain a classification of maximal proper subgroups $H$ of
$\PGSp(4, \mathbb{F}_q)$. It is similar to
the one above, except that cases 7) and 8) change according to the
relation between $H$ and $G$, given by the exact sequence:
$$ 1 \rightarrow  G \rightarrow  H \rightarrow  \{  \pm 1 \} \rightarrow 1 .$$

\subsection{The image of inertia at $\ell$}
In Weissauer's construction, the representations $\rho_{f,\lambda}$
are realized in the 'tale cohomology of the Siegel variety $X_1 $
of level $1$. Applying results of Faltings and Chai-Faltings (see [6],
[2]),
 Urban proves  (see [20]) a  result that combined with Fontaine-Laffaille
  theory (see [7]) gives:
\begin{prop}
\label{teo:Urban} Let $f$ be a Siegel cusp form as in theorem
\ref{teo:Weissauer}. Then $\rho_{f,\lambda}$ is crystalline and has
Hodge-Tate weights $\{ 2k-3 , k-1 , k-2 , 0 \}$. Moreover if
$\ell -1 > 2k - 3  $, then we have the following possibilities for the
action of the inertia group at $\ell$: $\bar{\rho}_{f,\lambda}|_{I_\ell} \simeq $
$$
 \pmatrix{
  1 & * & * & * \cr
  0 & \chi^{k-2} & * & * \cr
  0 & 0 & \chi^{k-1} & * \cr
  0 & 0 & 0 & \chi^{2k-3} \cr}  , \;
  \pmatrix{
  \psi^{2k-3} & 0 & * & * \cr
  0 & \psi^{(2k-3)\ell} & * & * \cr
  0 & 0 & \psi^{(k-2)+(k-1)\ell} & 0 \cr
  0 & 0 & 0 & \psi^{(k-1)+(k-2)\ell} \cr}  ,
  $$
   $$
  \pmatrix{
  \psi^{2k-3} & 0 & * & * \cr
  0 & \psi^{(2k-3)\ell} & * & * \cr
  0 & 0 & \chi^{k-2} & * \cr
  0 & 0 & 0 & \chi^{k-1} \cr} , \;
  \pmatrix{
  1 & * & * & * \cr
  0 & \chi^{2k-3} & * & * \cr
  0 & 0 & \psi^{(k-2)+(k-1)\ell} & 0 \cr
  0 & 0 & 0 & \psi^{(k-1)+(k-2)\ell} \cr}  ,
  $$
where $\psi$ denotes a fundamental character of level $2$.

\end{prop}

\section{Study of the images}
\subsection{Reducible case: $1$-dimensional constituent}
Let $f$ be a Siegel cusp form verifying all the conditions of theorem
\ref{teo:Weissauer}.
Suppose that for a prime $\lambda $ in $ E$
the representation $\bar{\rho}_{f,\lambda}$ is reducible with a
$1$-dimensional sub(or quotient) representation. The representation
being unramified outside $\ell$ we conclude from proposition
\ref{teo:Urban} (assume $\ell > 2k-2$) that this $1$-dimensional constituent is:
$  \chi^i$, $ i = 0, k-2, k-1$ or $2k-3$.\\
 As in [5],
we will use the fact that symplectic representations have ``reciprocal"
roots, i.e., that the roots of the characteristic polynomial of
$\rho_{f,\lambda}(\Frob \; p)$ come in pairs $ \alpha , p^{2k-3}/
\alpha$.\\
{\bf Remark:} We have also used this fact in proposition \ref{teo:Urban} to
discard many cases, in particular the cases involving fundamental characters of
level $3$ or $4$.\\
Therefore we can assume that  $\chi^i$ is a root of
the characteristic polynomial $Pol_p(x)$ of $\bar{\rho}_{f, \lambda}( \Frob \; p)$
$$ Pol_p(x) = x^4 - a_p x^3 + (a_p^2 - a_{p^2} - p^{2k-4} ) x^2 - a_p p^{2k-3} x +
p^{4k-6} ,$$
 for every $p \neq \ell$, with $i= 0$ or $k-1$.

 \medskip
 \noindent For $i=0$ we obtain the following congruence:
$$ b_p - a_p (p^{2k-3}+1) + p^{4k-6} + 1 \equiv 0 \pmod{\lambda}
 \quad \qquad (4.1)
 $$
for every $p \neq \ell$,
 where $b_p = a_p^2 - a_{p^2} - p^{2k-4}$ is the quadratic coefficient
of $Pol_p(x)$.\\
Weissauer proved (see [21]) that any genus $2$ Siegel cusp form that is
not in the image of the Saito-Kurokawa lift (and more generally any irreducible
 cuspidal
automorphic representations $\pi$ of $\GSp(4 , \mathbb{A})$  such that
$\pi $ is not Cuspidal Associated to Parabolic and $\pi_\infty$ is a discrete
series representation), i.e., not a Maass spezialform, verifies the
generalized Ramanujan conjecture, so  the roots of $Pol_p(x)$ have all
the same absolute value $\sqrt{p^{2k-3}}$.\\
From now on we will only work with Siegel cusp forms that are not
Maass spezialformen, with the twofold intention of obtaining large
images and using Ramanujan's conjecture.
 This is the case for the examples computed by
Skoruppa (cf. [15], [16]), where the  Siegel cusp forms that are not
Maass spezialformen are called
``interesting" (as we will call them).\\
Therefore, if $f$ is an interesting Siegel cusp form, using the bounds
for the absolute values of $a_p$ and $b_p$ that follow from
Ramanujan's conjecture we conclude that congruence (4.1) is not an
equality (for  large enough $p$). Therefore, it can only hold for
finitely many primes $\ell$.

\medskip
\noindent Similarly, for $i=k-1$ we obtain the congruence:
$$ p^{2k-4} (1+p^2)- a_p p^{k-2}(1+p)+b_p \equiv
0 \pmod{\lambda}
 \quad \qquad (4.2)
 $$
for every $p \neq \ell$.\\
Again, Ramanujan's conjecture implies that this is not an equality for
large enough $p$. Thus, the reducible case with $1$-dimensional
constituent can only hold for finitely many primes.

\subsection{Reducible case: related $2$-dimensional constituents}

Suppose that, after semi-simplification, $\bar{\rho}_{f,\lambda}$ decomposes
as the sum of two $2$-dimensional irreducible Galois representations:
$ \bar{\rho}_{f,\lambda} \cong \pi_1 \oplus \pi_2$. Assume also that these two
constituents are related, i.e., if
$\alpha , \beta$ are the roots of the characteristic polynomial of
$\pi_1 (\Frob \; p)$, then $p^{2k-3}/\alpha , p^{2k-3}/ \beta$
 are the roots of that of
$\pi_2 (\Frob \; p)$. If not, then : $\alpha =
p^{2k-3}/ \beta$, so  $\det(\pi_1) = \det(\pi_2) = \chi^{2k-3}$; this case will
 be studied in the next subsection.\\
 We apply proposition \ref{teo:Urban} (assume $\ell > 2k-2$).
 There are two  cases to consider:
 \begin{itemize}
 \item
 Case 1 : $\det(\pi_1) = \chi^{k-1} , \quad
  \det(\pi_2) = \chi^{ 3k-5} $.
  \item
Case 2: $ \det(\pi_1) = \chi^{k-2} , \quad
   \det(\pi_2) =  \chi^{3k-4} $.
   \end{itemize}
$\bullet$ Case 1: In this case we have the factorization:
$$ Pol_p(x) \equiv (x^2- A x + p^{3k-5} ) \;
( x^2 - \frac{A x}{ p^{k-2} } + p^{k-1}) \quad \pmod{\lambda} .$$
 Eliminating
$A $ from the equation, we obtain:
$$ (b_p - p^{k-1} - p^{3k-5}) \; ( p^{k-2}+1)^2 - a_p^2 \; p^{k-2}
\equiv 0
 \pmod{\lambda} \quad \qquad (4.3)$$
 for every $p \neq \ell$.\\
From the bounds on the coefficients that follow from Ramanujan's conjecture,
 we see that for large enough
$p$ this is not an equality. Thus, only finitely many $\ell$
can satisfy (4.3).\\
$\bullet$ Case 2: This case is quite similar to the previous one. We
start with:
$$ Pol_p(x) \equiv (x^2- A x + p^{3k-4} ) \;
( x^2 - \frac{A x}{  p^{k-1} } + p^{k-2}) \quad \pmod{\lambda} .$$
From this:
$$ (b_p - p^{k-2} - p^{3k-4}) \; ( p^{k-1}+1)^2 - a_p^2 \; p^{k-1}
\equiv 0
 \pmod{\lambda} \quad \qquad (4.4)$$
 for every $p \neq \ell$.\\
From the bounds on the coefficients we conclude that for an
interesting Siegel cusp form the reducible case with related
two-dimensional constituents can only hold for finitely many primes.

\subsection{The remaining reducible case with Serre's conjecture}

As explained above, in the remaining reducible case we have:
 $ \bar{\rho}_{f,\lambda}^{ss} \cong \pi_1 \oplus \pi_2$ with $\det(\pi_1)=
 \det(\pi_2)=\chi^{2k-3}$. Assume $\ell > 2k-2$.
 In proposition \ref{teo:Urban}
 we have given a description of
 $\bar{\rho}_{f,\lambda}|_{I_\ell}$. This gives for $\pi_1|_{I_\ell}$,
 $\pi_2|_{I_\ell}$:
$$ \pmatrix{
  1 & *  \cr
  0 &  \chi^{2k-3} \cr } \quad   \mbox{or} \quad
\pmatrix{
  \psi^{2k-3} & 0 \cr
  0 & \psi^{(2k-3)\ell} \cr} $$
  and
  $$ \pmatrix{
  \chi^{k-2} & *  \cr
  0 &  \chi^{k-1} \cr } \quad   \mbox{or} \quad
\pmatrix{
  \psi^{(k-2)+(k-1)\ell} & 0 \cr
  0 & \psi^{(k-1)+(k-2)\ell} \cr} ,$$
  respectively.
Besides, both two-dimensional representations are unramified outside
$\ell$.
At this point we invoke Serre's conjecture ($3.2.4_?$) (see [14])
that gives us a control on $\pi_1$ and $\pi'_2 := \chi^{-k+2} \otimes
\pi_2$.
 Both representations
should be modular of weights $2k-2$ and $2$, respectively, and level $1$; i.e.,
 there exist two cusp forms
$f_1 , f_2$ with:
$$ \bar{\rho}_{f_1,\lambda} \cong \pi_1 , \;  \bar{\rho}_{f_2,\lambda} \cong \pi'_2 ,
\quad f_1 \in S_{2k-2}(1) , \; f_2 \in S_2(1) $$
 (we are assuming that $\pi_1$ and $ \pi_2$ are
  irreducible; otherwise the results of section 4.1
  apply). Both cusp forms have trivial nebentypus.\\
But $S_2(1)=0$ and we obtain a contradiction.\\
We conclude, assuming
Serre's conjecture, that the reducible case with unrelated two
dimensional constituents cannot happen if $\ell > 2k-2$ for genus two Siegel
cusp forms for the full Siegel modular group $\Sp(4, \mathbb{Z})$ verifying the
conditions of theorem \ref{teo:Weissauer}.\\
{\bf Remark:} In all reducible cases (sections 4.1 , 4.2 and 4.3) we have
considered reducibility over $\bar{\mathbb{F}}_{\lambda}$.

\subsection{Groups with a reducible index $2$ subgroup}
If $G_\lambda$, the image of $\bar{\rho}_{f,\lambda}$, corresponds to an
irreducible subgroup inside (its projective image) some of the maximal subgroups in
 cases 3), 4) and 5) of Mitchell's classification,
 there is a normal subgroup of index $2$ of
 $G_\lambda$:
$$ 1 \rightarrow M_\lambda \rightarrow G_\lambda
 \rightarrow \{ \pm 1 \} \rightarrow 1 ,$$
where the subgroup $M_\lambda$ is reducible (not necessarily over
  $\mathbb{F}_\lambda $).\\
For $\ell> 2k-2$ we apply the description of
 $\bar{\rho}_{f,\lambda}|_{I_\ell}$ given in proposition \ref{teo:Urban}.
  Observe that $\chi^{k-2}$ has order larger than $2$ because $\ell - 1 >
 2k-3$, and the order of
 $\psi^{(k-2)+(k-1)\ell}$ is clearly a multiple
 of $\ell +1$, therefore larger than $2$.  Then we conclude that the image of $I_\ell$
  is contained in $M_\lambda$. \\
Therefore, if we take the quotient $G_\lambda /M_\lambda$
we obtain a representation
$$ G_\mathbb{Q} \rightarrow C_2 $$
 whose kernel
is a quadratic field unramified everywhere;  we thus obtain a contradiction. \\
We conclude that the image of $\bar{\rho}_{f,\lambda}$ never falls in
this case if $\ell > 2k-2$.

\subsection{The stabilizer of a twisted cubic}
In this case all upper-triangular matrices are of the
 form:
 $$  \pmatrix{
  a^3 & * & * & * \cr
  0 & a^2 d & * & * \cr
  0 & 0 & a d^2 & * \cr
  0 & 0 & 0 & d^3 \cr} .$$
Let us compare  this with the four possibilities for the image of the
 inertia subgroup at $\ell$ for $\ell > 2k-2$ given
 in proposition \ref{teo:Urban}.\\
 In the first case we see that this inertia
subgroup has the required
form only if: $\chi^{2k-3} \equiv \chi^{3k-6}$ , $\chi^{2k-3} \equiv \chi^{ 3k-3} $,
$ \chi^{3k-4} \equiv 1 $ or $ \chi^{3k-4} \equiv \chi^{6k-9}$. But this implies
$\ell-1 \mid k-3$ , $ \ell -1 \mid k$, $ \ell-1 \mid 3k-5$ or
$\ell-1 \mid 3k-4$, respectively. With the restriction $\ell > 2k-2$,
we see that this is impossible.\\
In the second case, we see that the inertia subgroup has the required
form only if:
$$ \psi^{(3k-4)+(3k-5)\ell} \equiv  \psi^{6k-9} \; \mbox{or} \;
  \psi^{(6k-9)\ell} ,$$
   or
    $$
  \psi^{(2k-3)+(4k-6)\ell} \equiv  \psi^{(3k-3)+(3k-6)\ell} \; \mbox{or} \;
  \psi^{(3k-6)+(3k-3)\ell}  .$$
   But this implies
  $\ell+1 \mid 3k-5$, $ \ell +1 \mid 3k-4$, $ \ell+1 \mid k$ or
$\ell+1 \mid k-3$, respectively. Again, this is impossible for $\ell >
2k-2$.\\
In the third case, we see that the inertia subgroup has the required
form only if:
$$ \psi^{ (2k-3 )(1+2\ell)} \equiv \chi^{3k-6} \; \mbox{or} \;
\chi^{3k-3} ,$$
 or
 $$ \chi^{3k-4} \equiv \psi^{6k-9} \; \mbox{or} \; \psi^{(6k-9)\ell} .$$
 The first congruence implies $\ell +1 \mid 2k-3$, which is impossible
for $\ell > 2k-2$. The second congruence implies:
$$ \ell^2 -1 \mid (-3k+5)+(3k-4)\ell , \qquad \quad (\ddagger) $$
which is impossible for every $\ell$ because $\ell -1 \mid (-3k+4) + (3k-4)
\ell$.\\
In the fourth and last case, if the inertia subgroup has the required form then:
$$ \psi^{(3k-4)+(3k-5)\ell} \equiv 1 \; \mbox{or} \; \chi^{6k-9} ,$$
 or
$$ \chi^{2k-3} \equiv \psi^{ (3k-6) + (3k-3) \ell  } \; \mbox{or} \;
   \psi^{ (3k-3) + (3k-6) \ell  } .$$
    The first congruence
   again gives ( $\ddagger$).\\ The second implies $ \ell +1 \mid (3k-6)+ (3k-3)\ell$,
   and from this: $\ell +1 \mid 3$.\\
 We conclude that the image does
not fall inside the stabilizer of a twisted cubic for any $\ell >
2k-2$.\\
{\bf Remark:} We are working with even weight. For the case of weight $3$,
the image of inertia may fall inside the stabilizer of a twisted
cubic.

\subsection{The exceptional groups}
We call  exceptional groups those appearing in cases 7) and 8) of the classification.
In these cases, comparing the exceptional group
 $H \subseteq \PGSp(4, \mathbb{F}_\lambda)$ or
$G \subseteq \PSp(4, \mathbb{F}_\lambda)$, its order and structure,
 with the fact that the image of $\mathbb{P}(\bar{\rho}_{f,\lambda})$
 contains the image of $ \mathbb{P}( \bar{\rho}_{f,\lambda}|_{I_\ell})$,
 and applying proposition \ref{teo:Urban}  (assuming $\ell > 2k-2$) we
 conclude that these cases can never happen for any $\ell >2k-2$. If $k < 6$
 we have to demand also that $\ell > 7$, but there is no interesting
 Siegel cusp forms for the full Siegel modular group
 $\Sp(4, \mathbb{Z})$ of weight smaller than $20$.

\subsection{Smaller symplectic groups}

To obtain the maximal possible image we need to impose the following
condition:
\begin{defin}
\label{teo:untwisted} \rm We say that a genus two Siegel cusp form $f$ is
untwisted if there exists a prime $p$ such that $a_p \neq 0$ and
$\mathbb{Q}(b_p) = E $, where $b_p = a_p^2- a_{p^2} - p^{ 2k-4}$.
\end{defin}
Alternatively, we can ask  the number field $E$ to be generated by
$ a_p^2$ for some prime $p \;$ (*),
 both conditions are enough to exclude for almost every prime
 the case where the image
is in a smaller symplectic group. This second condition is essentially the
same as that required in the case of classical modular forms to exclude the
case of inner twists (cf. [12], [4]).\\
Suppose that the Siegel cusp form $f$ is untwisted and let $p$ be a
prime as in definition \ref{teo:untwisted}. Let $\ell$ be a prime
different from $p$ such that
$$\ell \nmid \Norm(a_p) \quad \mbox{and}  \quad \ell \nmid \disc(b_p) . \qquad \; (4.5) $$
 Let $\lambda $ be a prime of $E$ lying above
$\ell$. We will show that for such a $\lambda$ the projective image of
$\bar{\rho}_{f,\lambda}$ cannot fall in a proper symplectic subgroup
of $\PGSp(4, \mathbb{F}_\lambda)$, i.e., it cannot be a subgroup as in
cases 9) and 10) of the classification.\\
Assume  the contrary. Then, there exists $c \in
\mathbb{F}_\lambda^*$ such that the polynomial
$$ x^4 -  c  a_p x^3 + c^2 b_p x^2 - c^3 a_p p^{2k-3} x + c^4 p ^{4k-6}
 \quad \mod \; \lambda$$
 has its coefficients in $ \mathbb{F}_{\lambda'}$ with
 $ [\mathbb{F}_\lambda : \mathbb{F}_{\lambda'}] > 1$. \\
 This implies that $ \frac {c^3 a_p} { c  a_p} = c^2 \in
  \mathbb{F}_{\lambda'}^* $ (recall that $\ell \nmid \Norm(a_p)$).
Thus, $ c^2 b_p / c^2 = b_p \in  \mathbb{F}_{\lambda'}$, contradicting
the assumption that $\mathbb{Q}(b_p) = E $ and $\ell \nmid \disc(b_p)$.
This contradiction proves that the image does not fall in a proper
symplectic group, if $\ell $ verifies (4.5).\\
{\bf Remark:} In a similar way, one can prove under condition (*) that the
image does not fall in a proper symplectic group, for every $\ell$
such that $\ell \nmid \Norm(a_p)$ and $ \ell \nmid \disc(a_p^2)$.\\

\subsection{Unconditional results without Serre's conjecture}
Without using Serre's conjecture  we can show that the images are ``as
large as possible"
for large density sets of primes, applying tricks similar to those used
in [5]. For $\lambda$ a prime in $E$, let $r$ be its
residue class degree:
 $\mathbb{F}_{\ell^r} = \mathbb{F}_\lambda$. Serre's conjecture was used
  to eliminate the
following two possibilities:\\

 i) The image of $\bar{\rho}_{f,\lambda}$ is contained in
 $$  \{  A \times B \in \GL(2,\mathbb{F}_\lambda) \times
  \GL(2,\mathbb{F}_\lambda) :
   \det(A)=\det(B)= \chi^{2k-3} \} . \qquad (4.6) $$
 ii) The image  of $\bar{\rho}_{f,\lambda}$ is contained in
 $$  \{ M \in \GL(2, \mathbb{F}_{\ell^{2r}} ) : \det(M)= \chi^{2k-3} \}
 \quad (4.7) $$

\medskip
\noindent The inclusion of this group in $\GSp(4, \mathbb{F}_\lambda)$ is given by
 the map:\\
 $M \rightarrow \diag(M, M^{\Frob})$, where $\Frob$ is the non-trivial
 element in $\Gal( \mathbb{F}_{\ell^{2r}} / \mathbb{F}_{\ell^r})$.

 \medskip
\noindent $\bullet$ Case i):
  We will use the standard factorization (see [15]):
  $$  Pol_p(x) = (x^2- (a_p/2 + \sqrt{d_p}) x + p^{2k-3})
  (x^2- (a_p/2 - \sqrt{d_p}) x + p^{2k-3}) ,$$
  where $d_p = -3/4 a_p^2 + a_{p^2} + p^{2k-4} + p^{2k-3}$.\\
  We impose the condition:
  $$ \mathbb{Q}(d_p) = E \quad \mbox{and} \quad  \sqrt{d_p} \notin
  E  . \quad \quad (4.8) $$
 Then, if we restrict to the primes $\lambda $ in $ E$ such that
 $$
  \ell \nmid \disc(d_p) \quad \mbox{and} \quad
   d_p \notin
 (\mathbb{F_\lambda})^2 \qquad \quad (4.9) $$
  case i) cannot hold, because the matrices
   $A$ and $B$ in (4.6) would have their traces in $\mathbb{F}_{\ell^{2r}}
    \smallsetminus
   \mathbb{F}_\lambda$.\\
   $\bullet$ Case ii): Assume $\ell > 2k-2$ and apply proposition
   \ref{teo:Urban}. Because $M= M^{\Frob}$ for any $M \in \GL(2,\mathbb{F}_\lambda) $, if the matrices in the image of inertia were in case ii) it should
hold:
$$ \{ 1, \chi^{2k-3}\} = \{ \chi^{k-2} , \chi^{k-1} \}  ,$$
$$      \{ \psi^{2k-3} , \psi^{(2k-3)\ell}  \} =
           \{   \psi^{(k-2)+(k-1)\ell}   ,  \psi^{(k-1)+(k-2)\ell}  \}   ,$$
$$      \{ \psi^{2k-3} , \psi^{(2k-3)\ell}  \} =
           \{   \chi^{k-2} , \chi^{k-1}  \}  \; \mbox{or}$$
           $$      \{ 1 , \chi^{2k-3}  \} =
           \{   \psi^{(k-2)+(k-1)\ell}   ,  \psi^{(k-1)+(k-2)\ell}  \}  .$$
The first of these equalities
implies $\ell - 1 \mid k-2$ or $ \ell -1 \mid k-1$.
The second implies $ \ell + 1 \mid k-2$ or $ \ell +1 \mid k-1 $.
The third implies $\ell + 1 \mid 2k-3$, and the fourth implies
$\ell +1 \mid (k-2) + (k-1)\ell$, impossible
for every $\ell$.
We conclude that no $\ell >
2k-2$ can verify any of these equalities, so  case ii) can never happen
for such an $\ell$.

\subsection{Conclusion}

Having gone through all cases of the classification (see section 3.1),
 we conclude that
under certain conditions on a Siegel cusp form $f$ the images of the
attached Galois representations are ``as large as possible" for almost
every prime, to be more precise:
\begin{teo}
\label{teo:mainresult} Let $f$ be
  a Siegel cusp form
 (and Hecke eigenform)  of even weight $k$
 for the full Siegel modular group $\Sp(4,\mathbb{Z})$ verifying the
 conditions:\\
$\bullet$ Multiplicity one.\\
 $\bullet$ Not a Maass spezialform. \\
 $\bullet$ Untwisted (see definition \ref{teo:untwisted}).\\
Assume Serre's conjecture $(3.2.4_?)$. Then the images of the Galois
representations
$\rho_{f,\lambda}$ are $A_\lambda^k$, for almost every $\lambda$.\\
Without assuming Serre's conjecture, take a prime $p$ such that (4.8)
holds. Then for all but finitely many primes verifying (4.9) the
 image of $\rho_{f,\lambda}$ is $A_\lambda^k$.
\end{teo}
{\bf Remark:} We are using the fact (see section 2) that the maximality of the images of
$\rho_{f,\lambda}$,
$\bar{\rho}_{f,\lambda}$ and
$\mathbb{P}( \bar{\rho}_{f,\lambda} )$ are equivalent. In section 6 we will
prove, without assuming Serre's conjecture,
 that the above result applies to almost every prime.\\
The result, with or without the assumption of Serre's
conjecture, is effective, and we can compute the exceptional primes
given a Siegel cusp form verifying the above conditions.\\
In fact, assuming Serre's conjecture, all exceptional primes (more exactly,
 a finite
 set
containing all exceptional primes) are
computed using equations (4.1), (4.2), (4.3), (4.4) and (4.5) (this last condition
implies in particular that $\ell$ is unramified in $E/ \mathbb{Q}$) applied to
 more than one $Pol_p(x)$ with the restrictions $\ell > 2k-2 $
.\\
Without assuming Serre's conjecture, we also compute several $d_p$
such that (4.8) holds and
find the finitely many exceptional primes among the primes verifying
(4.9) for any of these $d_p$.\\

\section{Examples}

We start by applying our method to an example of Skoruppa (see [15],
[16]), the example already investigated in [3]. This example
consists
of a Siegel Hecke eigenform $ \Upsilon$ of weight $k=28$ for the full Siegel
modular group $\Sp(4,\mathbb{Z})$ such that the space spanned by the
Galois conjugates of $\Upsilon$ is the complement of the space spanned
by the Eisenstein-Klingen series and the Maass spezialformen, so
the automorphic representation corresponding to $\Upsilon$ has
multiplicity one.\\
 Using the results of
[15] and the tables in [16], we can compute the characteristic
polynomials $Pol_p(x)$ for $p=2,3$ and $5$ as  in [3].\\
The field $E$ generated by the eigenvalues of $\Upsilon$ is the cubic field
generated by some root $\alpha$ of $x^3-x^2-294086x-59412960$.\\
The required Fourier coefficients of $\Upsilon$:
$$ a(1,1,1); a(2,2,2); a(3,3,3); a(4,4,4); a(5,5,5); a(1,1,7); a(3,3,7); a(1,1,19)$$
 are given in [16] in terms of
$\alpha$.
From these coefficients, using the formulas in [15],  the first eigenvalues of
$\Upsilon$: $a_2 , a_4 , a_3 , a_9 , a_5 , a_{25}$;
 are computed (see [3]).\\
  These values determine all the coefficients of the characteristic
  polynomials $Pol_p(x)$ for $p=2,3$ and $5$. From now on let us  impose
  $\ell > 2k-2 = 54$.\\
  We compute the exceptional primes falling in a reducible case
  using formulas (4.1) to (4.4) with $p=2$ and $3$.
   To do so just take the norm in these
  equations and discard the primes $\ell$ dividing the greatest
  common divisor of the  norms obtained for $p = 2$ and $3$. In the
  example we found no exceptional primes $\ell > 54$ at this step.\\
  As a matter of fact, in this and other steps,
  we have to exclude the
  ``primes dividing denominators" from consideration. These are the
  primes dividing the denominators of the norm of $a_2$, $a_3$, $a_5$,
  $a_4$, $a_9$ or $a_{25}$. In the example, they form the set:
  $$ \mathbb{D}= \{ 2,3, 17 , 2063 , 8841304187 , 1646767084367711  \} .$$
We use equation (4.5) for $p=2,3$ to bound the set of exceptional
primes falling in a smaller symplectic group ($a_2$ and $a_3$ are not zero and
$\mathbb{Q}(b_2) = \mathbb{Q}(b_3) = E $). Excluding the primes in
$\mathbb{D}$, we know that any such exceptional prime $\ell$ has to
divide the greatest common divisor of the two numbers:
$$ \numerator(\Norm(a_p)) \cdot \numerator(\disc(b_{p})) , \quad p=2,3 .$$
As the only possible exceptional primes we obtain the primes that ramify in
$E$:
$$ \mathbb{R} = \{  5,13,73693, 1418741 \} .$$
Therefore, we have the following result:
\begin{teo}
\label{teo:el28Serre} Let $\Upsilon$ be the  Siegel
cusp form of weight $28$ computed by Skoruppa. It is untwisted and has
multiplicity one. Let $\rho_{\Upsilon , \lambda}$
be the Galois representations attached to $\Upsilon$ constructed by
Weissauer. Then, if we assume Serre's conjecture $(3.2.4_?)$,
it follows that for every $\ell>53$, $\ell \notin \mathbb{D} \cup \mathbb{R}$
and every $\lambda $ in $E$ above $\ell$
the image of  $\rho_{\Upsilon , \lambda}$ is $A_\lambda^{28}$. In
particular, except for the finitely many primes excluded, the groups
$\PSp(4, \mathbb{F}_{\ell^r})$, if $r$ is even, and $\PGSp(4, \mathbb{F}_{\ell^r}
)$,
if $r$ is odd, (where $r$ denotes the residual degree of $\lambda$)
are realized as Galois group over $\mathbb{Q}$, and the
corresponding Galois extension only ramifies at $\ell$.
\end{teo}
$ $\\
Without assuming Serre's conjecture, we want to obtain an effective
unconditional result for the images of the $\rho_{\Upsilon,\lambda}$.
With this aim, we compute the values $d_p$, for $p=2,3$ and $5$, and
check that (4.8) holds for the three of them. Therefore, applying (4.9)
three times, we obtain:

\begin{teo}
\label{teo:el28sinSerre} Let $\Upsilon$,  $\rho_{\Upsilon , \lambda}$
be as in theorem \ref{teo:el28Serre}.
Then, for every $\ell>53$, $\ell \notin \mathbb{D} \cup \mathbb{R}$
and every $\lambda \in E$ above $\ell$
such that:
$$ \ell \nmid \disc(d_p) \quad \mbox{and } d_p \notin ( \mathbb{F}_\lambda)^2 , \quad
p=2,3, \; \mbox{or} \; 5 ,$$
the image of  $\rho_{\Upsilon , \lambda}$ is $A_\lambda^{28}$ (with the
same consequences for inverse Galois theory  as in theorem
\ref{teo:el28Serre}).
\end{teo}
\noindent
For the primes $\ell$ inert in $E/ \mathbb{Q}$
it is easy to see that if
$$ \ell \nmid \disc(d_p) \quad \mbox{and} \quad
 \Bigl(  \frac {\Norm(d_p)} {\ell} \Bigr) = -1  $$
  for $p =2,3$ or $5$, then
 $  d_p \notin ( \mathbb{F}_\lambda)^2  $. This gives an easier criterion
 for these primes.\\
 The first  primes $\ell > 53$ inert in $E$
are:
$$ 59, 67 , 71 , 101 , 103 , 137 , 151 ,157 , 181 , 191 , 197 .$$
We apply the above result to them and we verify that except for $151$,
the other ten primes in this list verify the required conditions for
some $p \leq 5$. Therefore we obtain the following corollary:

\begin{cor}
\label{teo:los10} The groups $ \PGSp(4 , \mathbb{F}_{\ell^3})$ are
Galois groups over $\mathbb{Q}$ for $\ell =  59 , 67 , 71 ,101
, 103 , 137  , 157 , 181 , 191 , 197$.
\end{cor}

\subsection{The first interesting Siegel cusp form}

Now we apply the algorithm to $\Upsilon20$, the first example
of interesting Siegel cusp form for the full Siegel modular group.
 $\Upsilon
20$ has weight $k=20$ and rational Fourier coefficients. Its
first eigenvalues were computed by Skoruppa in
[15]. The dimension of the space of interesting Siegel cusp forms of
 weight $20$ is $1$, therefore $\Upsilon 20$ has multiplicity
 one.\\
 We use the values of the $a_p$, $a_{p^2}$ for $p \leq 7$.\\
Applying equation (4.1) with $p=2$ and $3$ we find no exceptional
primes
at this step.\\
Applying equation (4.2) with $p=2,3$ and $5$ we obtain as (possible)
exceptional primes $\ell =2,3,5,7,11$. According to the results of
Skoruppa (cf. [15]) these primes are in fact exceptional because
 $\Upsilon 20$ is congruent modulo these $\ell$ (and only them) to
a Maass spezialform, so equation (4.2) holds, for these $\ell$,
 for every $p$. \\
Applying equation (4.3) with $p=2,3,5$ and $7$, we obtain as (possible)
exceptional primes $\ell = 2,3,5,29,71$.\\
Recall that with our algorithm we can show non-exceptionality only
for primes  $\ell > 2k-2 = 38$. Here we have the first prime $\ell >37$
in this example seeming to be exceptional: $\ell = 71$.  Kurokawa
proved in [9] that $\Upsilon 20$ is congruent modulo $71$ to an
Eisenstein-Klingen series, so $71$ is in fact exceptional and
(4.3) holds, for $\ell = 71$, for every $p$.\\
Applying equation (4.4) with $p=2,3,5$ only gives $\ell =2,3,5$ as
possible exceptional primes.

\medskip
\noindent This concludes the application of the algorithm to this example, proving
in particular that the image of the Galois representations $\rho_\ell$
attached to $\Upsilon 20$ are ``as large as possible" for every
$\ell > 71$, under the assumption of Serre's conjecture. To obtain an
unconditional result, we use the first values of $d_p$:\\
$d_2= 2^{14} \cdot 3^2 \cdot
 7 \cdot 13 \cdot 19  \cdot 241 $,  $\;
d_3= 2^6 \cdot 3^{10} \cdot
19 \cdot 47 \cdot 150628997 $, \\ $
d_5= 2^8 \cdot 3^2 \cdot 5^6 \cdot
 19 \cdot 47 \cdot 1396135808326877 $, \\ $
d_7=
  2^8 \cdot 3^6 \cdot 7^6 \cdot
29 \cdot 1097 \cdot 41713094306662453 $.
\begin{teo}
\label{teo:el20}  Let $\Upsilon 20$ be the  Siegel
cusp form of weight $20$ computed by Skoruppa. It  has
multiplicity one and $E = \mathbb{Q}$. Let $\rho_{ \ell}$
be the Galois representations attached to $\Upsilon 20$ constructed by
Weissauer. Then, if we assume Serre's conjecture $(3.2.4_?)$
it follows that for every $\ell>37$, $\ell \neq 71$
the image of  $\rho_{ \ell}$ is $A_\ell^{20}$. In
particular for every $\ell > 37$, $\ell \neq 71$,  the groups
 $\PGSp(4, \mathbb{F}_{\ell} )$
are realized as Galois group over $\mathbb{Q}$, and the
corresponding Galois extension of $\mathbb{Q}$ only ramifies at $\ell$.\\
The same conclusion holds unconditionally (without assuming
Serre's conjecture) if we also impose $ ( \frac{ d_p}{\ell }) = -1 $ for
$p=2,3,5$ or $7$.
\end{teo}
{\bf Remark:} In particular, we have proved that the images
are
``as large as possible" for an explicitly given
  set of primes of Dirichlet density $\frac{15 }{16
} \;$.

\section{Improving theorem \ref{teo:mainresult}}
\subsection{Pseudo-representations}

We reproduce Wiles' definition of
pseudo-representations
and some of its properties (cf. [23], [19]). This is a tool used to
construct odd two-dimensional representations.\\
Let $A$ be a ring, $G$ a group, and $c \in G$ an element of order two.
\begin{defin}
\label{teo:pseudoW} \rm A pseudo-representation of $G$ defined over $A$ is a
quadruple $\tau=(a,d,t,x); \; a,d,t : G \rightarrow A, \; x: G \times G
 \rightarrow A $ satisfying the conditions:\\
 {\rm (i)} $ 2 a_{g \cdot g'} = a_g a_{g'} + x_{g,g'} , \;
 2 d_{g \cdot g'} = d_g d_{g'} + x_{g',g}$\\
 {\rm (ii) } $ a_g = t_g + t_{cg} , \; d_g = t_g - t_{cg}$\\
 {\rm (iii)}  $t_1 = 2, \; t_c = 0, \; x_{g,c} = x_{c,g} = 0$\\
 {\rm (iv)  } $ x_{g,g'} x_{h,h'} = x_{g,h'} x_{h,g'} , \\
  4 x_{gh,g'h'} = a_g a_{h'} x_{h,g'} + a_{h'} d_h x_{g,g'} + a_g d_{g'} x_{h,h'}
   + d_h d_{g'} x_{g,h'}$.\\
\end{defin}
Define the trace of $\tau$ by $\trace(\tau) = t$, and the determinant
by \\ $\det(\tau)(g)= a_g d_g - x_{g,g}$.\\
{\bf Remarks:} (i) If $G$ and $A$ have a topology, all maps are assumed to be
continuous.\\
(ii) If $\rho: G \rightarrow \GL(2,A), \; \rho= (\alpha , \beta , \gamma , \delta)$
 is a representation defined over $A$, with $\rho(c)= (1,0,0,-1)$
  (an odd representation)
 it gives rise to a pseudo-representation $\tau = (a,d,t,x)$, where
 $a=2 \alpha, \; d= 2 \delta, \; t= \alpha+\delta, \; x_{g,g'} = 4
 \beta_g \gamma_{g'} $.\\
 (iii) A pseudo-representation is determined by its trace. If $A' \subset A$
  is a subring, if $\tau$ is defined over $A$ and $\trace(\tau)$ takes
  values in $A'$, then $\tau$ is defined over $A'$.\\
The following are two basic properties that we will use in the
sequel (cf. [23], [19]):
\begin{lema}
\label{teo:patching} (``patching lemma") Let $(U_i)_{i>0}$ be a
sequence of ideals in $A$ such that $\bigcap_{i} U_i = \{ 0\}$. For
each $i$, let $\tau_i$ be a pseudo-representation defined over $A/U_i$.
Assume that there exist a dense subset $\Sigma$ of $G$ and
a function $T: \Sigma \rightarrow A$ such that for any $i$,
$T \; \mod \; U_i  = \trace(\tau_i)$. Then, there exists a
pseudo-representation $\tau$ defined over $A$ such that for any $i$,
$\tau \; \mod \; U_i = \tau_i$.
\end{lema}

\begin{teo}
\label{teo:Wiles} If $A$ is a field of characteristic different from
$2$, any pseudo-representation $\tau$ corresponds to an odd two dimensional
 representation
$\rho$.
\end{teo}

\subsection{From residual to $\lambda$-adic reducibility}
Let $f$ be a genus two Siegel cusp form verifying the conditions of
theorem \ref{teo:mainresult}, and consider the Galois representations
$\rho_\lambda := \rho_{f,\lambda}$. Suppose that for an infinite set of primes
$ \{ \lambda_i \}_{i=1}^\infty = \Lambda$ of $E$ the residual
representations $\bar{\rho}_{\lambda_i}$ are reducible
with irreducible two dimensional constituents of the same determinant.\\
 Taking a subset of
$\Lambda$ if necessary, we assume that the density of the set $\mathbb{L}$
of rational primes $\ell_i$ such that $\lambda_i \mid \ell_i$ is $0$.\\
We want to prove that with this hypothesis of residual reducibility the
whole family of $\lambda$-adic representations $\rho_\lambda$ has to be
reducible  (and a fortiori residually reducible). We
assume for the moment
that the decompositions:
$$\bar{\rho}_{\lambda_i}^{ss} = \sigma_i \oplus \pi_i $$
induce a decomposition on the set $\{ a_p \}$ of traces of the family
$\rho_\lambda$: $ a_p = u_p + v_p \; ( \bigodot )$ for every
 $p \notin \mathbb{L}$ that
is independent of $\lambda_i$. That is to say, there exists such a
decomposition (fixed) with $\{ u_p \; \mod \; \lambda_i \}$ equal to the
set of traces of $\sigma_i$ (restricted to $p \notin \mathbb{L}$) and
$\{ v_p \; \mod \; \lambda_i \}$ equal to the
set of traces of $\pi_i$ (restricted to $p \notin \mathbb{L}$),
for every $\lambda_i \in \Lambda$.\\
Fix a prime $\lambda$ of $E$. As usual, let  $\ell$ be the rational prime
 such that $\lambda \mid \ell$. \\
For any set $S$ of rational primes let us denote by $\mathbb{Q}^S$ the
maximal extension of $\mathbb{Q}$ unramified outside $S$ and
$G^S = \Gal( \mathbb{Q}^S / \mathbb{Q})$. We will work only
with Galois representations unramified outside $ \mathbb{L}+\ell : =
 \mathbb{L} \bigcup \{ \ell \} $, therefore we  take the quotient
 $G^{\mathbb{L}+\ell}$ of $G_\mathbb{Q}$ and work with representations
 of this quotient group.  \\
Consider the subgroup $W_\lambda$
of $   G^{\mathbb{L}+\ell}$ algebraically generated by the elements
$ \{ \Frob \; p\}_{ p \notin \mathbb{L} + \ell}$.\\
 Let us consider the function $F: W_\lambda \rightarrow \overline{\mathbb{Q}} $
 defined on the generators of $W_\lambda$ by: $F(\Frob \;p) = u_p$.\\
 $ $\\
 Recall that the traces $a_p$ of $\rho_\lambda$ are the eigenvalues of
 the Siegel cusp form $f$. They
  generate the number field $E$. The decomposition of the
  characteristic polynomials $Pol_p$ introduced in section 4.8 shows
  that it must hold: $u_p = a_p / 2 \pm \sqrt{d_p}$ ,  so  we have:
  $$ E' := \mathbb{Q}( \{ u_p \}) = E ( \{ \sqrt{d_p} \}) .$$
Moreover, from the definition of $d_p$ (see section 4.8) and standard
properties of the eigenvalues of $f$ it is known that
all $a_p$ and $a_{p^2}$ belong to $\mathcal{O}_E [ \frac{1}{N} ]$
where $N$ is a finite product of primes and
all $u_p$ belongs to $\mathcal{O}' :=\mathcal{O}_{E'}[ \frac{1}{2N} ]$, where
$\mathcal{O}_{X}$ denotes  ``ring of algebraic integers of $X$".\\
$ $\\
Therefore, the image of $F$ is inside $\mathcal{O}'$:
$$ F: W_\lambda \rightarrow \mathcal{O}' . \quad \qquad (6.1)$$
For each $\lambda_i \in \Lambda$ we take a prime $\lambda'_i$ of $E'$
above $\lambda_i$. Composing the function $F$ with the reduction modulo
$\lambda'_i$ for every $\lambda_i \in \Lambda$ gives the functions:
$$ F_{\lambda'_i}: W_\lambda \rightarrow \mathbb{F}_{\lambda'_i}$$
with $F_{\lambda'_i}(\Frob \; p) = u_p \; \mod \; \lambda'_i$, for every
$p \notin \mathbb{L} + \ell$. This is
the restriction to $W_\lambda$ of the pseudo-representation $ \tilde{\sigma}_i$
associated to $\sigma_i$ , therefore $F_{\lambda'_i}$ is a
$2$-dimensional pseudo-representation of $W_\lambda$
 for every $\lambda_i \in \Lambda$.\\
 {\bf Remark:} For every $\lambda_i$, $\sigma_i$ gives a representation of
 $G^{\mathbb{L}+\ell}$ because it is unramified outside $\ell_i $.\\
Now we apply the patching lemma
 (see lemma \ref{teo:patching}) and we conclude that $F$ as in (6.1) is
a $2$-dimensional pseudo-representation. We have already seen (see for example
section 4.3) that the representations $\sigma_i$ have to be odd.
 An application of theorem \ref{teo:Wiles} proves that there exists a
representation
$$ \sigma: W_\lambda \rightarrow \GL(2, E')$$
with $\trace \; \sigma = F$, i.e.,
$$ \trace ( \sigma (\Frob \; p) ) = u_p$$
for every $p \notin \mathbb{L}+ \ell$.\\
The same argument shows that by patching the representations $\pi_i$ we
construct a  representation:
$$ \pi : W_\lambda \rightarrow \GL(2, E')$$
with
$$ \trace ( \pi (\Frob \; p) ) = v_p $$
for every $p \notin \mathbb{L}+ \ell$.\\
Now let us compare the semisimplification
 $({\rho_\lambda}|_{W_\lambda})^{ss}$ with
$\sigma \oplus \pi$, as representations of $W_\lambda$
with image in $\GSp(4, \overline{E})$.\\
{\bf Remark:} $\rho_\lambda$ gives a representation of
 $G^{\mathbb{L}+\ell}$ and $W_\lambda$ because it is unramified outside $\ell $.\\
 Both representations
 have the
same trace at every element of $W_\lambda$,
so they are conjugated. We conclude that ${\rho_\lambda}|_{W_\lambda}$ is
reducible for every $\lambda$.\\
But
the Cebotarev density theorem implies that
the subgroup $W_\lambda$ is dense in $G^{\mathbb{L}+\ell}$, because
$\mathbb{L}$ has density $0$. Therefore, being continuous,
the representation $\rho_\lambda$ has to be reducible for every prime $\lambda$
of $E$, of
the form
 $$ \rho_{f,\lambda}^{ss} = \sigma_\lambda
 \oplus \pi_\lambda $$
with $\det ( \sigma_\lambda ) =
 \det ( \pi_\lambda ) = \chi^{2k-3}$.

 \medskip
 \noindent It remains to remove the hypothesis $ ( \bigodot )$. To do
 this, one should consider the tensor product representations
$ \sigma_i \otimes \pi_i$ for every $\lambda_i \in \Lambda$. These four
dimensional irreducible representations can be patched as in the above construction
 because each of
them has trace at $\Frob \; p $ equal to $\{ u_p v_p \; \mod \;
\lambda_i \}$ (here we need Taylor's version of the theory of pseudo-representations
 and in particular of lemma \ref{teo:patching} valid for the four
  dimensional case, see [17]). The resulting characteristic $0$
  representation will be tensor-decomposable: $\sigma \otimes \pi$, and
  using the uniqueness of tensor decompositions in the irreducible
  case, one concludes that for the two dimensional representations
  $\sigma$ and $\pi$ it holds:
  $$ \trace(\sigma(\Frob \; p)) + \trace(\pi(\Frob \; p)) = u_p + v_p  $$
 and the proof continues as before. Thus, the following result holds:

\begin{teo}
\label{teo:lareducibilidad} Let $f$ be a genus $2$ Siegel cusp form
verifying the conditions of theorem \ref{teo:mainresult}. Let $\rho_{f,\lambda}$
be the corresponding family of Galois representations.
 Suppose that the residual representations
$\bar{\rho}_{f,\lambda_i}$ are reducible for infinitely many primes $\lambda_i$
of $E$. Then, for every prime $\lambda$ the $\lambda$-adic
representation $\rho_{f,\lambda}$ is reducible.
\end{teo}
Let $f$ be a Siegel cusp form verifying the conditions of theorem
\ref{teo:mainresult}. Assume that we are in the situation of theorem
\ref{teo:lareducibilidad}. Then we have: $\rho_{f,\lambda}^{ss} = \sigma_\lambda
 \oplus \pi_\lambda$ for every $\lambda$, with $\det ( \sigma_\lambda ) =
 \det ( \pi_\lambda ) = \chi^{2k-3}$. Suppose that condition (4.8) is
 verified by some prime $p$. As was shown in section 4.8, this implies
 that for infinitely many $\lambda$ with
 $$ \ell \nmid \disc(d_p) \quad \mbox{and} \quad d_p \notin ( \mathbb{F}_\lambda)^2$$
the representation $\bar{\rho}_{f,\lambda}$ cannot be reducible.
Thus we obtain a contradiction with theorem \ref{teo:lareducibilidad}.\\
Combining this with the results of  section 4 we conclude:
\begin{teo}
\label{teo:super}
Let $f$ be a genus $2$ Siegel cusp form and Hecke eigenform of even
weight $k$ for the Siegel modular group $\Sp(4, \mathbb{Z})$ with
multiplicity one, untwisted and not  a Maass spezialform.\\
Suppose that there is a prime $p$ with $\mathbb{Q}(d_p)= E$ and
$\sqrt{d_p} \notin E$. Then the image of the Galois representation
$ \rho_{f,\lambda}$ is $A_\lambda^k$ for almost every $\lambda$.
\end{teo}
Remark: This result applies in particular to the examples considered in the
previous section.
%
\section{Bibliography}\
\noindent [1] {\it H. Blichfeldt},  Finite collineation groups, (1917)
University of Chicago Press

\noindent [2] {\it C. Chai, G. Faltings}, Degeneration of Abelian Varieties, (1990)
Springer-Verlag

\noindent [3] {\it M. Dettweiler, U. Khn, S. Reiter},  On Galois
representations via Siegel modular forms of genus $2$, (2000) preprint

\noindent [4] {\it L. Dieulefait, N. Vila},  Projective linear groups as Galois
 groups over $\mathbb{Q}$ via modular representations, J. Symbolic
 Computation {\bf 30} (2000) 799-810

\noindent [5] {\it L. Dieulefait}, Explicit determination of the images of the
 Galois representations attached to abelian surfaces with $\End(A) =
  \mathbb{Z}$, preprint (2001)

\noindent [6] {\it G. Faltings},  Crystalline cohomology and Galois
representations, in Algebraic Analysis, Geometry and Number Theory,
Proceedings of JAMI Inaugural Conference, (1990) John Hopkins Univ. Press

\noindent [7] {\it  J.M. Fontaine, G. Laffaille},  Construction de repr'sentations
$p$-adiques, Ann. Scient. c. Norm. Sup., $4^e$ s'rie, t. {\bf 15} (1982)
547-608

\noindent [8] {\it J. Hirschfeld}, Finite projective spaces of three
dimensions, (1985) Clarendon Press - Oxford

\noindent [9] {\it N. Kurokawa}, Congruences between Siegel Modular Forms of
Degree Two, Proc. Japan Acad. {\bf 55}, Ser. A (1979) 417-422

\noindent [10] {\it H. Mitchell}, The subgroups of the quaternary abelian
linear group, Trans. Amer. Math. Soc. {\bf 15} (1914) 379-396

\noindent [11] {\it T. Ostrom}, Collineation groups whose order is prime to the
characteristic, Math. Z. {\bf 156} (1977) 59-71

\noindent [12] {\it K. Ribet}, On $\ell$-adic representations attached
 to modular
forms II, Glasgow Math. J. {\bf 27} (1985) 185-194

\noindent [13] {\it J-P. Serre}, Oeuvres, vol. 4, 1-55, (2000) Springer-Verlag

\noindent [14] {\it J-P. Serre}, Sur les repr'sentations modulaires de degr'
$2$ de $\Gal(\bar{\mathbb{Q}} / \mathbb{Q})$, Duke Math. J. {\bf 54}
(1987) 179-230

\noindent [15] {\it N. Skoruppa}, Computations of Siegel modular forms of genus
two, Math. Computation {\bf 58} (1992) 381-398

\noindent [16] {\it N. Skoruppa}, Siegel modular eigenforms of even weight on
the full Siegel modular group of genus $2$, available at:
http://thor.math.u-bordeaux.fr/ $\sim$modi

\noindent [17] {\it R. Taylor}, Galois representations associated to Siegel
modular forms of low weight, Duke Math. J. {\bf 63} (1991) 281-332

\noindent [18] {\it R. Taylor}, On the $\ell$-adic cohomology of Siegel
threefolds, Invent. Math. {\bf 114} (1993) 289-310

\noindent [19] {\it J. Tilouine}, Galois representations congruent to those
coming from Shimura Varieties, Proc. Sympos.
Pure Math. {\bf 55} Part 2 (1994) 625-638

\noindent [20] {\it E. Urban}, Selmer groups and the Eisenstein-Klingen Ideal,
 preprint (1998)

\noindent [21] {\it R. Weissauer}, The Ramanujan Conjecture for genus two
Siegel modular forms (An application of the Trace Formula), preprint
(1994)

\noindent [22] {\it R. Weissauer}, Four dimensional Galois representations,
preprint (2000)

\noindent [23] {\it A. Wiles}, On nearly ordinary $\lambda$-adic representation
 associated to modular forms, Invent. Math. {\bf 94} (1988) 529-573


\end{document}